\documentstyle[12pt, epsf, epsfig, amssymb, amstext, amstex,amsthm, amsfonts]{article}
\input xy
\input amssymb.sty
\xyoption{all}
\begin{document}

\newtheorem{thm}{Theorem}[section]
\newtheorem{prop}[thm]{Proposition}
\newtheorem{cor}[thm]{Corollary}
\newtheorem{lem}[thm]{Lemma}
\newtheorem{conj}[thm]{Conjecture}
\newtheorem{exa}[thm]{Example}
\newtheorem{defn}[thm]{Definition}
\newtheorem{clm}[thm]{Claim}
\newtheorem{eex}[thm]{Exercise}
\newtheorem{obs}[thm]{Observation}
\newtheorem{note}[thm]{Notation}
\newtheorem{remark}[thm]{Remark}

\newcommand{\brem}{\begin{remark}}
\newcommand{\erem}{\end{remark}}
 
\newcommand{\ben}{\begin{enumerate}}
\newcommand{\een}{\end{enumerate}}
\newcommand{\blem}{\begin{lem}}
\newcommand{\elem}{\end{lem}}
\newcommand{\bcl}{\begin{clm}}
\newcommand{\ecl}{\end{clm}}
\newcommand{\bthm}{\begin{thm}}
\newcommand{\ethm}{\end{thm}}
\newcommand{\bpr}{\begin{prop}}
\newcommand{\epr}{\end{prop}}
\newcommand{\bco}{\begin{cor}}
\newcommand{\eco}{\end{cor}}
\newcommand{\bcon}{\begin{conj}}
\newcommand{\econ}{\end{conj}}
\newcommand{\bde}{\begin{defn}}
\newcommand{\ede}{\end{defn}}
\newcommand{\bex}{\begin{exa}}
\newcommand{\eexa}{\end{exa}}
\newcommand{\bexe}{\begin{exe}}
\newcommand{\eexe}{\end{exe}}
\newcommand{\bobs}{\begin{obs}}
\newcommand{\eobs}{\end{obs}}
\newcommand{\bnote}{\begin{note}}
\newcommand{\enote}{\end{note}}

\newcommand{\fg}{\Pi _1(D-K,u)}
\newcommand{\Z}{{\Bbb Z}}
\newcommand{\C}{{\Bbb C}}
\newcommand{\R}{{\Bbb R}}
\newcommand{\Q}{{\Bbb Q}}
\newcommand{\F}{{\Bbb F}}
\newcommand{\N}{{\Bbb N}}
\newcommand{\HE}{\overset{HE}{\backsim}}
\newcommand{\fnref}[1]{~(\ref{#1})}

\newenvironment{emphit}{\begin{itemize} \em}{\end{itemize}}

\begin{center}
\Large{\bf {Hurwitz Equivalence of Braid Group Factorizations Consisting of a Semi-Frame}}\\ 
\vspace{7mm}
\large{T. Ben-Itzhak and M. Teicher} \footnote{This work was partially supported by the Emmy Noether Research Institute for Mathematics (center of the Minerva Foundation of Germany), the Excellency Center "Group Theoretic Methods in the Study of Algebraic Varieties"  of the Israel Science Foundation, and by EAGER (EU network , HPRN-CT-2009-00099).\\}

\end{center}

%             ABSTRACT
%       --------------------

\vspace{15mm}
ABSTRACT. In this paper we prove certain Hurwitz equivalence properties in the braid group. Our main result is that every two factorizations of $\Delta _n ^2$ where the elements of the factorization are semi-frame are Hurwitz equivalent. The results of this paper are generalization of the results in \cite{B4}. We use a new presentation of the braid group, called the Birman-Ko-Lee presentation, to define the semi-frame structure. The main result of this paper can be applied to compute the BMT invariant of surfaces (presented in \cite{KuTe} or \cite{Te}). The BMT is the class of Hurwitz equivalent factorizations of the central element of the braid group. The BMT distinguish among diffeomorphic surfaces which are not deformation of each other.

%             Basic Definitions
%-------------------------------------

\section{Topological Background}

In this section we recall some basic definitions and statements from \cite{MoTe1}:\\
Let $D$ be a closed disk on $\R ^2$, $K \subset D$ finite set, $u \in \partial D$. Any diffeomorphism of $D$ which fixes $K$ and is the identity on $\partial D$ acts naturally on $\Pi _1 = \Pi _1(D-K,u)$. We say that two such diffeomorphisms of $D$ (which fix $K$ and equal identity on $\partial D$) are equivalent if they define the same automorphism on $\Pi _1(D-K,u)$. This equivalence relation is compatible with composition of diffeomorphism and thus the equivalence classes form a group.\\

% Bn[D,K] Definition. 
%---------------------

\bde \label{BraidGroupDef}
Braid Group $B_n [D,K]$.
\ede
\noindent Let $D,K$ be as above, $n = \# K$, and let $\mathcal{B}$ be the group of all diffeomorphisms $\beta$ of $D$ such that $\beta (K)=K$, $\beta |_{\partial D} = Id_{\partial D}$. For $\beta _1, \beta _2 \in \mathcal{B}$ we say that $\beta _1$ is equivalent to $\beta _2$ if $\beta _1$ and $\beta _2$ define the same automorphism of $\Pi _1(D-K,u)$. The quotient of $\mathcal{B}$ by this equivalence relation is called the braid group $B_n [D,K]$.\\
Equivalently, if we take the canonical homomorphism $\psi :{\mathcal{B}} \rightarrow {Aut}(\Pi _1(D-K,u))$, then $B_n [D,K] = Im(\psi )$. The elements of $B_n [D,K]$ are called braids.\\ 

\blem
If $K' \subset  D'$ in another pair as above with ${\# K'} = {\# K} = n$, then $B_n [D',K']$ is isomorphic to $B_n [D,K]$.\\
\elem

\noindent This gives rise to the definition of $B_n$:

\bde
$B_n = B_n[D,K]$ for some $D,K$ with $ \# K = n$.
\ede

\bnote
We use the notation $a[b]$ for conjugating $b$ to $a$.
\enote

% Half Twist
%------------
\bde \label{HalfTwistDef}
$H(\sigma)$, half-twist defined by $\sigma$.
\ede
\noindent Let $D$ and $K$ be defined as above. Let $a, b \in K$, $K_{a,b} = K \backslash \{ a,b \}$ and $\sigma$ be a simple path (without a self intersection) in $D \backslash \partial D$ connecting $a$ with $b$ such that $\sigma \bigcap K = \{ a,b \}$. Choose a small regular neighborhood  $U$ of $\sigma$ such that $K_{a,b} \bigcap U = \phi$, and an orientation preserving diffeomorphism $\psi : {\R}^2 \rightarrow {\C}^1$ such that
$\psi (\sigma) = [-1,1] = \{ z \in \C ^1 | Re(z) \in [-1,1], Im(z) = 0 \}$ and $\psi (U) = \{ z \in \C ^1||z|<2 \}$. Let $\alpha (r),r \geq 0$, be a real smooth monotone function such that $\alpha (r) = 1$ for $r \in [0,3/2]$ and $\alpha (r) = 0$ for $r \geq 2$. Define a diffeomorphism $h:\C ^1 \rightarrow \C ^1$ as follows: for $z \in \C ^1, z = r e^{i\phi} \quad \text{let} \quad h(z)= r e^{i(\phi +\alpha (r) \pi )}$. It is clear that the restriction of $h$ to $\{ z \in \C ^1 | |z| \leq 3/2 \}$ coincides with the $180 ^\circ$ positive rotation, and that the restriction to $\{ z \in \C ^1 | |z| \geq 2 \}$ is the identity map. The diffeomorphism $\psi ^{-1} \circ h \circ \psi$ induces a braid called half-twist and denoted by $H(\sigma)$\\

\blem \label{HtConj}
Let $H(\sigma _1)$ and $H(\sigma _2)$ be the half twists defined by the paths $\sigma _1$ and $\sigma _2$ respectively. Then $H(\sigma _1)[H(\sigma _2)] = H(H(\sigma _2)(\sigma _1))$.
\elem
\noindent
{\bf Proof:} III.1.0 of \cite{MoTe5}.

% Frame
%--------

\bde \label{frame}
Frame of $B_n[D,K]$
\ede
\noindent Let $K = \{ p_1,...,p_n \}$ and $\sigma _1,..., \sigma _{n-1}$ be a system of simple smooth paths in $D - \partial D$ such that $\sigma _i$ connects $p_i \text{ with } p_{i+1}$ and $L = \bigcup \sigma _i$ is a simple smooth path. The ordered system of half-twists $\{ H(\sigma _1 ),...,H(\sigma _{n-1}) \} $ is called a frame of $B_n [D,K]$

\bthm\label{relations} 
Let $\{ H_1,...,H_{n-1} \}$ be a frame of $B_n[D,K]$ then,  $B_n[D,K]$ is generated by $\{ H_i \} _{i=1} ^{n-1}$ and the following is a complete list of relations:
$$H_i H_j = H_j H_i \text{ if }|i-j|>1 \text{      (Commutative relation)}$$
$$H_i H_j H_i = H_j H_i H_j \text { if } |i-j|=1 \text{      (Triple relation)}$$
\ethm

\noindent Proof can be found for example in \cite{MoTe1}.\\
This theorem provides us with Artin's algebraic definition of the braid group. \\

% Band Generators
%------------------
\noindent
Recently, new presentation of $B_n$, called the band generators presentation has been introduced by Birman-Ko-Lee \cite{Birman}. The new presentation is based on Aritin's presentation. We give here the definition of the presentation slightly different by constructing the generators on a frame structure.\\

\begin{figure}[h]
\begin{center}
\epsfxsize=8cm
\epsfysize=4cm
\epsfbox{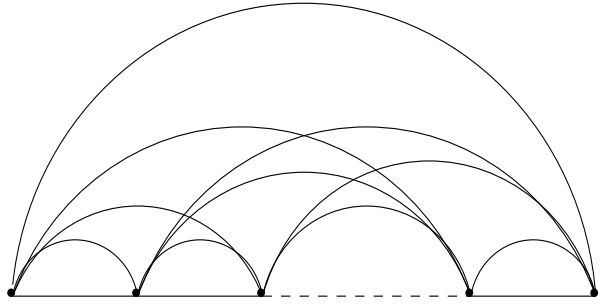}
\caption{Graph of Band Generators}
\label{FigBand}
\end{center}
\end{figure}

\bde \label{BandPresentation}
Band Generators Presentation
\ede
\noindent
Let $\{ H_1,...,H_{n-1} \}$ be a frame. The band generators corresponding to the frame  $\{ H_1,...,H_{n-1} \}$ are defined by the $(_2 ^n )$ generators: 
$$ a_{t,s} = (H_{t-1} H_{t-2} \cdots H_{s+1}) H_s ( H_{s+1} ^{-1} \cdots H_{t-2} ^{-1} H_{t-1} ^{-1}) \quad 1 \leq s < t \leq n $$
and the presentation's relations are:
\begin{equation} \label{BirmanRelation1}
a_{t,s} a_{s,r} = a_{t,r} a_{t,s} = a_{s,r} a_{t,r} \quad \text{ for all } t,s,r \text{ with } 1 \leq r < s < t \leq n
\end{equation}
\begin{equation} \label{BirmanRelation2}
a_{t,s} a_{r,q} = a_{r,q} a_{t,s} \quad \text{ if } (t-r)(t-q)(s-r)(s-q) > 0
\end{equation}
\noindent
The band generators contain the frame elements since $H_t = a_{t+1,t} $. An example of band generators corresponding to the trivial frame can be seen in Figure \ref{FigBand}. The paths in Figure \ref{FigBand} were computed using Lemma \ref{HtConj}.

\bde
Semi-Frame of $B_n[D,K]$
\ede
\noindent
We say that the system of half-twists $\{ H_i \} _{i=1} ^m$, is a semi-frame, if there exist a corresponding system of simple paths $\{ \sigma _i \} _{i=1} ^m$, such that, $\forall i$, $H_i = H(\sigma _i)$ and $\{ \sigma _i \} _{i=1} ^m$ satisfy the following conditions:
\begin{emphit} 
\item For each $\sigma _i$, $\sigma _i$ is a simple path in $D - \partial D$ connecting two points in $K$ and $| \sigma _i \bigcap K | = 2$.
\item There exist $n$ arcs $\{ \gamma _j \} _{j=1} ^n$, such that, $\gamma _j$ is a simple path in $ ( D - \partial D ) \bigcup \{ u \} $ connecting $u \in D$ with $p_j \in K$ and $\gamma _{j_0} \bigcap \gamma _{j_1} = \{ u \} \quad \forall j_0 \neq \j_1$.
\item Each $\gamma _j$ is intersecting $\bigcup _{i=1} ^m \sigma _i$ only at $ p_j $
\end{emphit}

\section{Definition of Hurwitz Moves}
%========================================================

% Hurwitz Move 

\bde
Hurwitz move on $G^m$ ($R_k, R_k ^{-1}$)
\ede
\noindent Let $G$ be a group, $\overrightarrow{t}=(t_1,...,t_m) \in G^m$. We say that $\overrightarrow{s}=(s_1,...,s_m) \in G^m$ is obtained from $\overrightarrow{t}$ by the Hurwitz move $R_k$ (or $\overrightarrow{t}$ is obtained from $\overrightarrow{s}$ by the Hurwitz move $R_k ^{-1}$) if 
$$ s_i = t_i \quad \text{ for } i \not= k, k+1,$$   
$$ s_k = t_k t_{k+1} t_k ^{-1}, \quad s_{k+1} = t_k.$$
 
% Factorization

\bnote
Factorization
\enote
\noindent Let $G$ be a group, a factorization $F$ of $g \in G$, is a list of factors $f_1,f_2,...,f_n$ ($f_i \in G$) such that the product $f_1 f_2...f_n$ is equal to $g$.\\
\noindent We use the notation $f_1 \cdot f_2 \cdots f_n$ for the factorization $f_1,f_2,...,f_n$ and $f_1 f_2...f_n$ for the product.

% Hurwitz Move on factorization.  

\bde
Hurwitz move on factorization 
\ede
\noindent Let $G$ be a group and $t \in G$. Let $t = t_1 \cdots t_m = s_1 \cdots s_m$ be two factorized expressions of $t$. We say that $s_1 \cdots s_m$ is obtained from $t_1 \cdots t_m$ by the Hurwitz move $R_k$ if $(s_1,...,s_m)$ is obtained from $(t_1,...,t_m)$ by the Hurwitz move $R_k$.\\

\bde
Hurwitz equivalence of factorization
\ede
\noindent The factorizations $s_1 \cdots s_m$, $t_1 \cdots t_m$\label{formin2} are Hurwitz equivalent if they are obtained from each other by a finite sequence of Hurwitz moves. The notation is  $t_1 \cdots t_m \overset{HE}{\backsim} s_1 \cdots s_m$.\label{formin1} \\

% Word
%----------

\bde
Word in $B_n$
\ede
\noindent
A word in $B_n$ is a representation of braid as a sequence of the generators and their inverses.

\bde
Positive Word in $B_n$
\ede
\noindent
Word in $B_n$ is positive if all generators are in positive powers.\\

% Positive Equal
%----------------

\section{Hurwitz Equivalence of Factorizations with Generators as Elements}
%=================================================================================================================
\noindent
The property proved in the following theorem is called the embedding property of a presentation, it was first proved by Garside in \cite{Garside} for the Artin presentation and recently was proved by Birman-Ko-Lee in \cite{Birman} for the band generators presentation.

\bthm \label{PositiveEqual}
Every two positive words (all generators with positive powers) which are equal, are transformable into each other through a finite sequence of positive words, such that each word of the sequence is positive and obtained from the preceding one by a direct application of one of the presentation relations.
\ethm
\noindent
{\bf Proof:} \cite{Garside},\cite{Birman}.\\

\noindent
In the following Lemma we show that the relations (\ref{BirmanRelation1}) and (\ref{BirmanRelation2}) also holds under the Hurwitz equivalence relation.

\blem \label{HurwitzRelations}
Let $\{ a_{t,s} \} \quad 1 \leq s < t \leq n-1 $ be a set of band generators as defined in \ref{BandPresentation}, then,\\
1. $a_{t,s} \cdot a_{s,r} \HE a_{t,r} \cdot a_{t,s} \HE a_{s,r}\cdot a_{t,r} \quad 
\text{ for all } t,s,r \text{ with } 1 \leq r < s < t \leq n$\\
2. $a_{t,s} \cdot a_{r,q} \HE a_{r,q} \cdot a_{t,s} \quad \text{ if } (t-r)(t-q)(s-r)(s-q) > 0$
\elem

\noindent
{\bf Proof:}\\
\noindent
1. $a_{t,s} \cdot a_{s,r} \HE a_{s,r}[a_{t,s}^{-1}] \cdot a_{t,s}$ by performing $R_1$. Since, $a_{t,s} a_{s,r} =(a_{s,r}[a_{t,s}^{-1}]) a_{t,s}$ as words in $B_n$, and  $a_{t,s} a_{s,r} = a_{t,r} a_{t,s}$ by relation (\ref{BirmanRelation1}) , we get that $a_{s,r}[a_{t,s}^{-1}] = a_{t,r}$ and therefore, $a_{t,s} \cdot a_{s,r} \HE a_{t,r} \cdot a_{t,s}$.\\
We use the same arguments to prove that $a_{t,s} \cdot a_{s,r} \HE a_{s,r}\cdot a_{t,r}$ by performing $R_1 ^{-1}$.\\
2. $a_{t,s} \cdot a_{r,q} \overset{R_1}{\rightarrow} a_{t,s} a_{r,q} a_{t,s} ^{-1} \cdot a_{t,s} =  a_{r,q} \cdot a_{t,s}$.\\

\noindent From Theorem \ref{PositiveEqual} and Lemma \ref{HurwitzRelations}, we get the following Theorem:

\bthm \label{sequence}
Let $\{ a_{t,s} \} \quad 1 \leq s < t \leq n-1 $ be a set of band generators of $B_n$ and $a_{i_1,k_1} \dots a_{i_p,k_p} = a_{j_1,l_1} \dots a_{j_p,l_p}$ two positive equal words (all $a_{i_r,k_r}$ and $a_{j_r,l_r}$ are band generators) then   $a_{i_1,k_1} \cdots a_{i_p,k_p}$ and $a_{j_1,l_1} \cdots a_{j_p,l_p}$ are Hurwitz equivalent as factorizations.
\ethm

\noindent 
{\bf Proof:} Applying Theorem \ref{PositiveEqual} on  $a_{i_1,k_1} \dots a_{i_p,k_p} = a_{j_1,l_1} \dots a_{j_p,l_p}$, we get a finite sequence of positive words $\{ W_r \} _{r=0} ^{q}$ s.t. $W_0 = a_{i_1,k_1} \dots a_{i_p,k_p}$, $W_q = a_{j_1,l_1} \dots a_{j_p,l_p}$ and $W_{r+1}$ is obtained from $W_r$ by a single application of one of the relations in \ref{BandPresentation}.\\
As we proved in Lemma \ref{HurwitzRelations}, each application of relation can be 'imitated' by a single Hurwitz move on the factorization. Thus, $W_0 \HE W_q$.

\bde
$\Delta _n ^2  \in B_n[D,K]$
\ede
$\Delta _n ^{2} = (H_1 \dots H_{n-1})^n$ where $ \{ H_i \} _{i=1} ^{n-1}$ is a frame.\\

\noindent We apply \ref{sequence} on $\Delta _n ^2$: 

\bco \label{sameFrameFactorization}
Let $\{ a_{t,s} \} \quad 1 \leq s < t \leq n-1 $ be a set of band generators corresponding to a frame  $ \{ H_i \} _{i=1} ^{n-1}$. Then all $\Delta _n ^2$ factorizations $a_{i_1, k_1} \cdots a_{i_{n(n-1)},k_{n(n-1)}}$ are Hurwitz equivalent to $(H_1 \cdots H_{n-1})^n$.\\
\eco

\blem \label{conjFrame1}
Let $\{ H({\sigma _i}) \} _{i=1} ^{n-1}$ be the frame which generates $B_n [D,K]$ and $H_{\zeta}$ a half twist. Then $\{ H({\sigma _i})[H_{\zeta}] \} _{i=1} ^{n-1}$ is also a frame.
\elem
\noindent {\bf Proof:} \cite{B4}.

\bthm \label{equivalenConj}
Let $\{ H_i \} _{i=1} ^{n-1}$ be a frame and let $b \in B_n[D,K]$, then the factorizations  $(H_1 \cdots H_{n-1})^n$ and $(H_1 [b] \cdots H_{n-1} [b])^n$ are Hurwitz equivalent.
\ethm
\noindent {\bf Proof:} \cite{B4}.

%
%		THE MAIN RESULT
%             =============================
%

\section{The Main Result}

\noindent In this section we proof the main theorem:
\bthm\label{MainTheorem}
Every two $\Delta ^2$ factorizations, $X_1 \cdots X_{n(n-1)}$ and $Y_1 \cdots Y_{n(n-1)}$ where $\{ X_i \} _{i=1} ^{n(n-1)}$ and $\{ Y_i \} _{i=1} ^{n(n-1)}$ are semi-frames are Hurwitz equivalent.
\ethm

\blem \label{MainLemma}
$\forall \gamma , \delta \in B_n$ any two $\Delta ^2$ factorizations, $a_{i_1,k_1}[\gamma ] \cdots a_{i_{n(n-1)},k_{n(n-1)}}[\gamma ]$ and $a_{j_1,l_1} [\delta ]\cdots a_{j_{n(n-1)},l_{n(n-1)}}[\delta ]$ where all $a_{i_r,k_r}$ and $a_{j_r,l_r}$ band generators, are Hurwitz equivalent.
\elem
\noindent
{\bf Proof:} Since, $\{ a_{t,s} \} $ are the band generators corresponding to the frame $\{ a_{t, t+1} \} _{t=1} ^ {n-1}$. By Lemma \ref{conjFrame1}, $\{ a_{t, t+1}[\gamma ] \} _{t=1} ^ {n-1}$ is also a frame and  $\{ a_{t,s}[\gamma ] \} $ are the corresponding band generators. Therefore, by Corollary \ref{sameFrameFactorization}:
\begin{equation} \label{rel1}
a_{i_1,k_1}[\gamma ] \cdots a_{i_{n(n-1)},k_{n(n-1)}}[\gamma ] \HE (a_{2,1}[\gamma ] \cdots a_{n,n-1}[\gamma ])^n
\end{equation}
\begin{equation} \label{rel2}
a_{j_1,l_1}[\delta ] \cdots a_{j_{n(n-1)},l_{n(n-1)}}[\delta ] \HE (a_{2,1}[\delta ] \cdots a_{n,n-1}[\delta ])^n
\end{equation}
By Theorem \ref{equivalenConj}, $(a_{2,1}[\gamma ] \cdots a_{n,n-1}[\gamma ])^n$ and $(a_{2,1}[\delta ] \cdots a_{n,n-1}[\delta ])^n$ are both Hurwitz equivalent to $(a_{2,1} \cdots a_{n,n-1})^n$ and therefore, equivalent to each other. From this and from relations (\ref{rel1}) and (\ref{rel2}) we get that:
$$ a_{i_1,k_1}[\gamma ] \cdots a_{i_{n(n-1)},k_{n(n-1)}}[\gamma ] \HE a_{j_1,l_1}[\delta ] \cdots a_{j_{n(n-1)},l_{n(n-1)}}[\delta ] $$

\noindent Lemma \ref{MainLemma} requires that the elements of the factorization are subset of the band generators corresponding to some frame. Since, conjugating $\gamma $ to a half twist is the same as operating $\gamma $ as a diffeomorphism on the path of the half twist (Lemma \ref{HtConj}), we are restricted to diffeomorphisms of subgraphs of the band generators in Figure \ref{FigBand}. The following Lemma completes the proof of the main theorem.

\begin{figure}[h]
\begin{center}
\epsfxsize=6cm
\epsfysize=6cm
\epsfbox{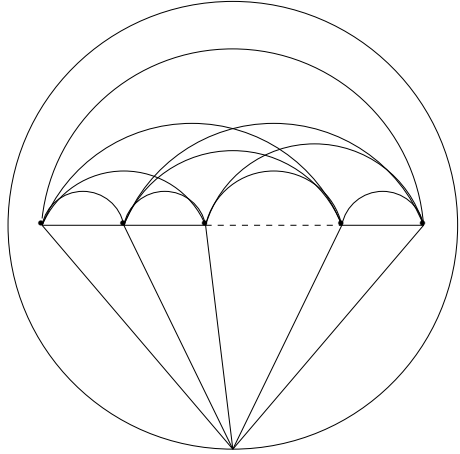}
\caption{The vertices can be connected by $n$ disjointed arcs to an external point.}
\label{FigArcs}
\end{center}
\end{figure}

\blem \label{GraphLemma}
A graph of generators is conjugated to a subgraph of the band generators graph if and only if there exist $n$ arcs which intersect with each other only at a point outside the graph, and intersect with the graph only at the vertices.
\elem
\noindent
{\bf Proof:} In the first direction, observe the graph of the band generators corresponding to the trivial frame. As shown in Figure \ref{FigArcs} the graph can be connected by $n$ disjointed arcs to an external point. Conjugating the generators in Figure \ref{FigArcs} to a braid is the same as conjugating to a list of generators. By Lemma \ref{HtConj} each conjugation is a diffeomorphism on the disk and therefore, we get that conjugating the braid to the generators is the same as operating diffeomorphism on the graph. Operating the same diffeomorphism on the arcs, we get $n$ disjointed arcs intersected with the graph only at the vertices.
The second direction of the Lemma was proved in \cite{Ko}.\\\\
Lemma \ref{GraphLemma} completes he proof of the main theorem.
$\hfill  \qed$\\

\brem
\erem
\noindent
If we consider factorizations where the elements of the factorization are not semi-frame, we will loose the embedding property, since it was proved in \cite{Ko} that the Artin presentation and the band generators presentation are the only presentations with the embedding property.

\end{document}